\documentclass[twoside, 11pt]{article}

\usepackage{mathrsfs,amsfonts,amsmath}
\RequirePackage{ifthen,calc}
\usepackage{theorem}
\usepackage[dvips]{color}

\def\E{{\mathbb E}}
\def\P{{\mathbb P}}

\def\N{{\mathbb R}^N}

 \catcode`@=11
 \def\@evenhead{\hbox to\textwidth{\footnotesize\rm\thepage \hfill
  {\it Convergence rates in precise asymptotics for
 a kind of complete moment convergence}}} % authors name

 \def\@oddhead{\hbox to \textwidth{\footnotesize{\it
 } \hfill\thepage}}% abbreviate title

 \renewcommand{\section}{\makeatletter
 \renewcommand{\@seccntformat}[1]{{\csname the##1\endcsname.}\hspace{0.45em}}
 \makeatother \@startsection
{section}%                                            the name
{1}%                                                  the level
{0pt}%                                                the indent
{\baselineskip}%                                      the beforeskip
{0.5\baselineskip}%                                   the afterskip
{\normalsize\bfseries\mathversion{bold}}}

\renewcommand{\subsection}{\makeatletter
 \renewcommand{\@seccntformat}[1]{{\csname the##1\endcsname.}\hspace{0.45em}}
 \makeatother \@startsection
{subsection}%                                            the name
{1}%                                                  the level
{0pt}%                                                the indent
{\baselineskip}%                                      the beforeskip
{0.5\baselineskip}%                                   the afterskip
{\normalsize\bfseries\mathversion{bold}}}
\catcode`@=12

\newtheorem{theorem}{\noindent Theorem}[section]
\newtheorem{lem}{\noindent Lemma}[section]

\newtheorem{prop}{\noindent Proposition}[section]

\theoremheaderfont{\normalfont\bfseries}
\theorembodyfont{\slshape} {\theorembodyfont{\rmfamily}} {\theorembodyfont{\rmfamily}
}
{\theorembodyfont{\rmfamily} } {\theorembodyfont{\rmfamily}
}
{\theorembodyfont{\rmfamily} } {\theorembodyfont{\rmfamily}
\newtheorem{rem}{\noindent Remark}[section]}
{\theorembodyfont{\rmfamily} } {\theorembodyfont{\rmfamily} } {\theorembodyfont{\rmfamily}
}
 \def\beqlb{\begin{eqnarray}}\def\eeqlb{\end{eqnarray}}
 \def\beqnn{\begin{eqnarray*}}\def\eeqnn{\end{eqnarray*}}
 \numberwithin{equation}{section}
\def\qed{\hfill$\square$\smallskip}
\def\N{{\mathbb{N}}}
\def\E{{\mathbb{E}}}
\def\P{{\mathbb{P}}}
\begin{document}
\title{\bf  Convergence rates in precise asymptotics for
 a kind of complete moment convergence
\footnotetext{\hspace{-5ex}
${[1]}$ School of Statistics, Shandong University of Finance and Economics,
Jinan 250014,  China. \\
 E-mail: kltgw277519@126.com (Kong, L.), mathdsh@gmail.com(Dai, H.: Corresponding author).
\newline
}}
\author{\small Lingtao Kong$^{1}$, Hongshuai Dai$^{1}$ }
\date{}

\maketitle

\begin{abstract}
In Liu and Lin (Statist. Probab. Letters, 2006), they introduced a kind of complete moment convergence
which includes complete convergence as a special case.  Inspired by the study of complete convergence,
in this paper,
we study the convergence rates of the precise asymptotics
for this kind of complete moment convergence and get  the corresponding convergence rates.
\end{abstract}

{\bf Keywords:} Convergence rate; precise asymptotics;
complete moment convergence
\vspace{2mm}

{\bf
MSC(2010): } 60F15 \quad 60G50
\vspace{2 mm}

\section{Introduction}

\setcounter{equation}{0}
Let $\{X, X_n, n\in\N\}$ be a sequence of i.i.d. random variables. Hsu and Robbins \cite{HR1947}  introduced the  complete convergence, and
proved that if $\E [X]=0$ and $\E [X^2]<\infty$, then
\begin{equation}\label{eq:1.1}
\sum_{n=1}^{\infty}
\P(|S_n|\geq\epsilon n)<\infty, ~\epsilon>0,
\end{equation}
where  $$S_n=\sum_{k=1}^n X_k, n\in \N.$$  Since then, it has attracted a lot of interest. For example,
Erd\H{o}s \cite{E1949, E1950} proved the necessity.
For more information, we refer to
Baum and Katz \cite{BK1965}, Davis \cite{D1968a, D1968b},
Lai \cite{L1974} and Gut \cite{G1980}.

Here we point out  that the sum in \eqref{eq:1.1} tends to infinity as $\epsilon\searrow0$. Due to this fact,
 it is interesting to study the  precise asymptotic problem, that is,
finding
an elementary function $f(\epsilon)$ such that $f(\epsilon)\sum_{n=1}^{\infty}
\P(|S_n|\geq\epsilon n)$ has a non-degenerate limit as $\epsilon\searrow0$.
 Many results about this topic have been established.
For example,  Heyde \cite{H1975} showed that
\begin{equation}\label{eq:1.2}
\lim_{\epsilon\searrow0}\epsilon^2\sum_{n=1}^{\infty}
\P(|S_n|\geq \epsilon n)=\E [X^2],
\end{equation}
whenever $\E [X]=0$ and $\E [X^2]<\infty$. For more results on precise asymptotic problems,
we refer to  Chen \cite{C1978},  Gut and Sp\u{a}taru \cite{GS2000a, GS2000b},
Li and Sp\u{a}taru \cite{LS2012}, Sp\u{a}taru \cite{S2004, S2014}
and the references therein.

Noting the Heyde's result (1.2),  it is natural  to
consider the rate of the convergence for the precise asymptotic
problem. This topic has been studied extensively.  For example,
 Klesov \cite{K1994}
 proved that  if $\E [X]=0$, $\E[X^2]>0$ and $\E[|X|^3]<\infty$, then
\begin{eqnarray}\label{eq:1.3}
\epsilon^2\sum_{n=1}^{\infty}\P(|S_n|\geq \epsilon n)-\E[X^2]=o(\epsilon^{\frac{1}{2}}).
\end{eqnarray}
Furthermore, by replacing  $\E[|X|^3]<\infty$ in  Klesov \cite{K1994} with $\E [|X|^q]<\infty$ for some $q\in(2,3]$,
He and Xie  \cite{HX2013}
got a much faster rate  for Heyde's result.
Recently, Gut and Steinebach \cite{GS2012,GS2013}, and Kong \cite{K2015} also extended the Klesov's \cite{K1994}
result, respectively, and got some new results.

As an extension of the complete convergence,  Liu and Lin \cite{LL2006}  introduced a new kind of complete moment convergence, and got
the precise asymptotic results for this kind of
complete moment convergence.
 These results read as follows.

\begin{theorem}\label{thm-1}
Let $\{X,X_n, n\in \N\}$ be a sequence of i.i.d.
random variables with partial sums $\{S_n,n\in\N\}$. For any $\epsilon>0$, set
\begin{eqnarray}
&&\lambda_1(\epsilon,p)=\sum_{n=1}^{\infty}\frac{1}{n^p}
\E\big[|S_n|^pI\{|S_n|\geq \epsilon n\}\big],~~\mbox{for any }p\in[0,2],\label{eq:1.4}\\
&&\lambda_2(\epsilon,\delta)=\sum_{n=2}^{\infty}\frac{(\log n)^{\delta-1}}{n^2}
\E \big[S_n^2 I\big\{|S_n|\geq\epsilon\sqrt{n\log n}\big\}\big],~~\mbox{for any }\delta\in(0,1].\label{eq:1.5}
\end{eqnarray}
(a)~For any $p\in[0,2)$, we have
\begin{equation}\label{eq:1.7}
\lim_{\epsilon\searrow0}\epsilon^{2-p}\lambda_1(\epsilon,p)=\frac{2\sigma^2}{2-p}
\end{equation}
if and only if  $\E [X]=0$ and $\E [X^2]=\sigma^2<\infty$.\\
(b)~For any $\delta\in (0,1]$, we have
\begin{eqnarray}\label{eq:1.8}
\lim_{\epsilon\searrow0}\epsilon^{2\delta}\lambda_2(\epsilon,\delta)=
\frac{\sigma^{2\delta+2}}{\delta}\E [|N|^{2\delta+2}]
\end{eqnarray}
if and only if
$\E [X]=0$, $\E [X^2]=\sigma^2$ and $\E [X^2(\log^{+}|X|)^{\delta}]<\infty$,
where $N$ denotes the  standard normal random variable and $\log^{+}x=\log (x\vee e)$.
\end{theorem}

Similar to the study of the rate of convergence in the complete convergence, it is also interesting to consider the convergence rates in precise
asymptotics for Liu and Lin's type complete moment convergence, which is the motivation of our work.  In this paper, we will carry out this work.

The rest of this paper is organized as follows.
In Section 2, we state the main results of this paper.
Section 3 is devoted to presenting the detailed proofs of these results.

We end this section with some notations.
We use $C$ to denote a positive constant whose value may vary from place to place.
Let $N$ denote the standard normal random variable and  $\Phi(x)=\P (|N|\geq x)$ for any $x\geq 0$.

\section{Main Results}
In this section, we state the main results of this paper. Before we do it,  we first introduce some notations.
Define
\begin{eqnarray}
B_{n,\theta}&=&\sum_{j=1}^nj^\theta-\frac{n^{\theta+1}}{\theta+1},
~~\mbox{ for any }-1<\theta<0;\label{eq:2.2}\\
C_{n,\delta}&=&\sum_{j=2}^n\frac{(\log j)^{\delta-1}}{j}-\frac{(\log n)^\delta}{\delta},
~~\mbox{ for any }0<\delta\leq 1.\label{eq:2.3}
\end{eqnarray}
We use $B_\theta$ and $C_\delta$ to denote the limits of the sequences
$\{B_{n,\theta},n\in \N \}$ and $\{C_{n,\delta},n\geq 2\}$, respectively.

\begin{rem}
Gut and Steinebach \cite{GS2013} proved the convergence
of the sequence $\{B_{n,\theta}; n\in \N\}$.
For the convergence
of the sequence
$\{C_{n,\delta}; n\geq 2\}$, we will give the proof in Section 3 below.
\end{rem}

Now, we state our main results. We have

\begin{theorem}\label{thm-2}
Let $\{X,X_n, n\in\N\}$ be a sequence of i.i.d.
random variables with partial sums $\{S_n,n\in\N\}$.
Suppose
that \begin{equation}\label{eq:2.2a}
\E [X]=0,~ \E [X^2]=\sigma^2>0,\mbox{ and }\E [|X|^q]<\infty \mbox{ for some }q\in(2,3].
\end{equation}

\begin{itemize}
\item[(a)] For any $p\in(0,2)$, we have
\begin{equation}\label{eq:2.5}
\lim_{\epsilon\searrow 0}\epsilon^{(\gamma-1)(2-p)}\Big[\epsilon^{2-p}\lambda_1(\epsilon,p)-\frac{2\sigma^2}{2-p}\Big]
=0,
\end{equation}
where $\gamma=\frac{q-p}{2q-2-p}$ and $\lambda_{1}(\epsilon,p)$ is given by \eqref{eq:1.4}.
\item[(b)] For any $0<\delta\leq 1$, we have
\begin{eqnarray}\label{eq:2.6}
\lim_{\epsilon\searrow0}(\log\frac{1}{\epsilon})^{-\delta}\epsilon^{-2\delta}\Big[\epsilon^{2\delta}\lambda_2(\epsilon,\delta)
-\frac{\sigma^{2\delta+2}}{\delta}\E[|N|^{2\delta+2}]\Big]=0,
\end{eqnarray}
where $\lambda_2(\epsilon,\delta)$  is defined by
\eqref{eq:1.5}.
\end{itemize}
\end{theorem}

\begin{rem}
Under the condition of Theorem \ref{thm-2}, Theorem \ref{thm-2} includes Theorem \ref{thm-1} as a special case.
\end{rem}

\begin{rem}
It is obvious that $0<\gamma<1$, since $p\in(0,2)$, $q\in(2,3]$.
\end{rem}

\begin{rem}
\eqref{eq:2.5} obtains the convergence rate  for $0<p<2$ in \eqref{eq:1.7}.
When $p=0$, then
$$\lambda_1(\epsilon,0)=\sum_{n=1}^{\infty}\frac{1}{n^2}\P (|S_n|\geq\epsilon n).$$
Under the assumptions of Theorem \ref{thm-2}, He and Xie \cite{HX2013} got that
\begin{eqnarray}\label{eq:2.4}
\epsilon^2\lambda_1(\epsilon,0)-\sigma^2=o(\epsilon^{q-2}).
\end{eqnarray}
Moreover,
Liu and Lin \cite{LL2006} also got the precise asymptotics for the case
$p=2$. However, the method we used in this paper seems unable to deal with the case $p=2$
and we will study this case in the future work.
\end{rem}

The proof of Theorem \ref{thm-2} is based on the following two propositions, one of which is concerned with
the Gaussian case (Proposition \ref{prop2.1}) and  the other is related to two Berry-Esseen type remainder terms
(Proposition \ref{prop2.2}).

\begin{prop}\label{prop2.1}
Suppose that  $\{X,X_n,n\in\N \}$ is a sequence of  i.i.d. normal random variables
with mean $0$ and variance $\sigma^2>0$. For any $-1<\theta<0$ and $0<\delta\leq1$,
$B_\theta$ and $C_\delta $ denote the limits of the sequences
$\{B_{n,\theta},n\in \N \}$ and $\{C_{n,\delta}; n\geq 2\}$, respectively.
\begin{itemize}
\item[(a)] For any $p\in (0,2)$, we have, as $\epsilon\searrow0$,
\begin{eqnarray}
\lambda_1(\epsilon,p)=\frac{2}{2-p}\epsilon^{p-2}\sigma^2+
B_{-\frac{p}{2}}\E[|N|^p]\sigma^p+O(\epsilon^p\log\frac{1}{\epsilon}).
\end{eqnarray}
\item[(b)] For any $\delta\in (0,1]$, we have, as $\epsilon\searrow0$,
\begin{eqnarray}
\lambda_{2}(\epsilon,\delta)
=\frac{1}{\delta}\E[|N|^{2\delta+2}]\epsilon^{-2\delta}\sigma^{2\delta+2}+C_{\delta}\E [N^2]\sigma^2+O(\epsilon^2).
\end{eqnarray}
\end{itemize}
\end{prop}

\begin{prop}\label{prop2.2}
Let $\{X,X_n;~n\in \N\}$ be a sequence of i.i.d. random
variables with $\E [X]=0$, $\E [X^2]=\sigma^2>0$, $\E [|X|^q]<\infty$ for some $q\in (2,3]$,
and partial sums $S_n=\sum_{k=1}^nX_k,~n\in \N$.
\begin{itemize}
\item[(a)] For any $p\in (0,2)$, we have,
\begin{eqnarray}\label{eq:2.8}
\lim_{\epsilon\searrow 0}\epsilon^{\gamma(2-p)}\sum_{n=1}^{\infty}\frac{1}{
n^{p}}\int_{\epsilon n}^{\infty}px^{p-1}\Big|\P(|S_n|\geq x)-\Phi\big( \frac{x}{\sqrt{n}}\big)\Big|dx=0,
\end{eqnarray}
where $\gamma=\frac{q-p}{2q-2-p}$.
\item[(b)] For any $\delta\in (0,1]$, we have,
 \begin{eqnarray}\label{eq:2.9}
\lim_{\epsilon\searrow 0}(\log\frac{1}{\epsilon})^{-\delta}
\sum_{n=2}^{\infty}\frac{(\log n)^{\delta-1}}{n^2}\int_{\epsilon\sqrt{n\log n}}^{\infty}
2x\Big|\P(|S_n|\geq x)-\Phi\big( \frac{x}{\sqrt{n}}\big)\Big|dx=0.
\end{eqnarray}
\end{itemize}
\end{prop}
\begin{rem} By Proposotions 4.2, 4.3, 5.2 and 5.3 in
Liu and Lin \cite{LL2006}, it's easy to  get the following results
$$\lim_{\epsilon\searrow 0}\epsilon^{2-p}\sum_{n=1}^{\infty}
\frac{1}{n^{p}}\int_{\epsilon n}^{\infty}px^{p-1}\Big|\P(|S_n|\geq x)-\Phi\big( \frac{x}{\sqrt{n}}\big)\Big|dx=0,$$
$$\lim_{\epsilon\searrow 0}\epsilon^{2\delta}\sum_{n=2}^{\infty}\frac{(\log n)^{\delta-1}}{n^2}\int_{\epsilon\sqrt{n\log n}}^{\infty}
2x\Big|\P(|S_n|\geq x)-\Phi\big( \frac{x}{\sqrt{n}}\big)\Big|dx=0.$$
We improve the above results by \eqref{eq:2.8} and \eqref{eq:2.9} with $0<\gamma<1$.
\end{rem}

\section{Proofs}
\setcounter{equation}{0}
In this section,  we give the detailed proofs of Theorem \ref{thm-2} and Propositions \ref{prop2.1} and \ref{prop2.2}.
 Without loss of generality,
we  assume $\sigma=1$ in this section.  In order to reach our aim, we first introduce some technical lemmas. The first one comes from Billingsley \cite{B1999}.

\begin{lem}\label{lem-3.1}
For $x$ large enough, We have
\begin{equation*}
\Phi(x)\sim \frac{2}{\sqrt{2\pi}x}e^{-\frac{x^2}{2}}.
\end{equation*}

\end{lem}

The following lemma comes from  Gut and Steinebach \cite{GS2013} which
shows the convergence of the sequence $\{B_{n,\theta}, n\in\N\}$ for any $\theta\in (-1,0)$.

\begin{lem}\label{lem-3.2}
We have, as $n\rightarrow\infty$,
\begin{equation}\label{eq:3.1}
B_{n,\theta}
=B_{\theta}+O(n^\theta),\mbox{~~~~~ for }-1<\theta<0,
  \end{equation}
  where $B_{n,\theta}$ is defined by \eqref{eq:2.2} and
  $B_{\theta}$ is a constant with $-\frac{1}{\theta+1}\leq B_\theta< \frac{\theta}{\theta+1}<0$.
\end{lem}

Similar to Lemma \ref{lem-3.2}, the following lemma shows the convergence of the sequence $\{C_{n,\delta}, n\geq2\}$ for any $\delta\in(0,1]$.

\begin{lem}\label{lem-3.3}
 For any $0<\delta\leq 1$, we have, as $n\rightarrow\infty$,
\begin{equation*}
C_{n,\delta}=C_{\delta}+O\Big(\frac{(\log n)^{\delta-1}}{n}\Big),
\end{equation*}
where $C_{n,\delta}$ is defined by \eqref{eq:2.3}
and $C_\delta$ is a constant with $-\frac{(\log 2)^{\delta}}{\delta}\leq C_{\delta}\leq 0$.
\end{lem}
{\it Proof:}
It follows from \eqref{eq:2.3} and  the mean value theorem that
\beqlb\label{eq:3.2}
C_{n+1,\delta}-C_{n,\delta}=&&
\frac{[\log (n+1)]^{\delta-1}}{n+1}-\frac{[\log (n+1)]^{\delta}-(\log n)^{\delta}}{\delta}\nonumber
\\
&&=\frac{[\log (n+1)]^{\delta-1}}{n+1}-\frac{(\log \xi_n)^{\delta-1}}{\xi_n}
\eeqlb
for some $\xi_n\in(n,n+1)$.

On the other hand, noting
\begin{equation}\label{eq:3.3}
C_{n,\delta}=\sum_{j=2}^{n}\int_{j-1}^{j}
\Big[\frac{(\log j)^{\delta-1}}{j}-\frac{(\log x)^{\delta-1}}{x}\Big]dx,
\end{equation}
 we get that
\beqlb\label{d-1}
0\geq C_{n,\delta}&\geq& \sum_{j=3}^{n}\Big[\frac{(\log j)^{\delta-1}}{j}-\frac{(\log (j-1))^{\delta-1}}{j-1}\Big]\nonumber
+\frac{(\log 2)^{\delta-1}}{2}-\int_{1}^2\frac{(\log x)^{\delta-1}}{x}dx\\
&=& \frac{(\log n)^{\delta-1}}{n}-\frac{(\log 2)^{\delta}}{\delta}\geq-\frac{(\log 2)^{\delta}}{\delta},
\eeqlb
since $f(x)=\frac{(\log x)^{\delta-1}}{x}$ is a
decreasing function.

Since the function $f(x)=\frac{(\log x)^{\delta-1}}{x}$ is
decreasing with $\delta\in (0,1]$, we obtain  that the sequence
$\{C_{n,\delta}, n\geq 2\}$ is also decreasing. By using this fact, we get from
the monotone bounded theorem  that $\lim_{n\rightarrow\infty}C_{n,\delta}$ exists,
 since \eqref{d-1} implies that $\{C_{n,\delta},n\geq 2\}$ is a bounded sequence.

For any $\delta\in(0,1]$,
let $C_{\delta}$  denote the limit of the sequence $\{C_{n,\delta},n\geq 2\}$. Moreover, from \eqref{d-1}
we get that
$-\frac{(\log 2)^{\delta}}{\delta}\leq C_{\delta}\leq 0$.

Given $m>n$, using the mean value theorem again, we have
\begin{eqnarray}
0>C_{m,\delta}-C_{n,\delta}
&=&\sum_{j=n+1}^m\Big[\frac{(\log j)^{\delta-1}}{j}-\frac{(\log j)^{\delta}-(\log (j-1))^{\delta}}{\delta}\Big]\nonumber\\
&=&\sum_{j=n+1}^m\Big[\frac{(\log j)^{\delta-1}}{j}-\frac{(\log \xi_j)^{\delta-1}}{\xi_j}\Big]
~~~~~~~~~\mbox{ for some }\xi_j\in(j-1,j)\nonumber\\
&>&\sum_{j=n+1}^m\Big[\frac{(\log j)^{\delta-1}}{j}-\frac{(\log (j-1))^{\delta-1}}{j-1}\Big]\nonumber\\
&=&\frac{(\log m)^{\delta-1}}{m}-\frac{(\log n)^{\delta-1}}{n}.\label{eq:3.4}
\end{eqnarray}
Letting $m\rightarrow\infty$ in \eqref{eq:3.4}, we
finish the proof of Lemma \ref{lem-3.3}.
\qed

Now we come to a point where we can prove the Propositions  \ref{prop2.1} and \ref{prop2.2}.  Next,  we first prove the Proposition \ref{prop2.1}.

\noindent{\it Proof of Proposition \ref{prop2.1}:}
Let  $\{X,X_n,n\in\N \}$ be a sequence of  i.i.d. normal random variables
with mean $0$ and variance $\sigma^2>0$. Since we assume $\sigma=1$,
we indeed deal with the standard normal random variable $N$.

We first prove the part (a) of Proposition \ref{prop2.1}.
 Note that, for $0<p<2$,
\begin{eqnarray}\label{eq:3.5}
\lambda_1(\epsilon,p)=\epsilon^p\sum_{n=1}^{\infty}\P(|S_n|\geq \epsilon n)
+\sum_{n=1}^{\infty}\frac{1}{n^p}\int_{\epsilon n}^{\infty}px^{p-1}\P(|S_ n|\geq x)dx.
\end{eqnarray}
On the other hand, from Klesov \cite{K1994}, we have
\begin{equation}\label{eq:3.6}
\sum_{n=1}^{\infty}\P(|S_n|\geq \epsilon n)=\frac{\sigma^2}{\epsilon^{2}}-\frac{1}{2}+o(1).
\end{equation}
Thus, in order to prove \eqref{eq:3.5}, we only need to consider the second term in \eqref{eq:3.5}.

The change of variable $t=\frac{x}{\sqrt{n}}$ yields
that\beqlb\label{d-2}
\frac{1}{n^p}\int_{\epsilon n}^{\infty}px^{p-1}\P(|S_ n|\geq x)dx\nonumber
&&=\frac{1}{n^{\frac{p}{2}}}
\int_{\epsilon\sqrt{n}}^{\infty}pt^{p-1}\Phi(t)dt\nonumber
\\&&=\frac{1}{n^{\frac{p}{2}}}
\sum_{j=n}^{\infty}\int_{\epsilon\sqrt{j}}^{\epsilon\sqrt{j+1}}pt^{p-1}\Phi(t)dt,
\eeqlb
since  $$\frac{S_n}{\sqrt{n}}\overset{d}{=} N$$ and
$$\Phi(t)=\P (|N|\geq t).$$
It follows from  the Fubini's theorem, \eqref{d-2} and Lemma \ref{lem-3.2} that
\begin{eqnarray}\label{eq:3.7}
&&\sum_{n=1}^{\infty}\frac{1}{n^p}\int_{\epsilon n}^{\infty}px^{p-1}\P(|S_ n|\geq x)dx\nonumber\\
&=&\sum_{j=1}^{\infty}\big(\sum_{n=1}^{j}\frac{1}{n^{\frac{p}{2}}}\big)
\int_{\epsilon\sqrt{j}}^{\epsilon\sqrt{j+1}}pt^{p-1}\Phi(t)dt\nonumber\\
&=&I_{31}+I_{32}+I_{33},
\end{eqnarray}
where
\beqnn
I_{31}&&=\frac{2p}{2-p}\sum_{j=1}^{\infty}j^{-\frac{p}{2}+1}
\int_{\epsilon\sqrt{j}}^{\epsilon\sqrt{j+1}}t^{p-1}\Phi(t)dt,
\\
I_{32}&&=B_{-\frac{p}{2}}
\sum_{j=1}^{\infty}\int_{\epsilon\sqrt{j}}^{\epsilon\sqrt{j+1}}pt^{p-1}\Phi(t)dt,
\\ I_{33}&&=\sum_{j=1}^{\infty}O(j^{-\frac{p}{2}})
\sum_{j=1}^{\infty}\int_{\epsilon\sqrt{j}}^{\epsilon\sqrt{j+1}}pt^{p-1}\Phi(t)dt.
\eeqnn

Now, we deal with $I_{31}$.
The integer mean theorem implies that for some $\xi\in(j,j+1)$
$$\int_{\epsilon\sqrt{j}}^{\epsilon\sqrt{j+1}}t^{p-1}\Phi(t)dt
=\epsilon^p\xi^{\frac{p-1}{2}}\Phi(\epsilon\sqrt{\xi})(\sqrt{j+1}-\sqrt{j}).$$
Furthermore, by the Taylor expansion, we have
\begin{eqnarray*}
\xi^{\frac{p-1}{2}}&=&j^{\frac{p-1}{2}}+O(j^{\frac{p-3}{2}}),\\
\Phi(\epsilon\sqrt{\xi})&=&\Phi(\epsilon\sqrt{j})+
\epsilon e^{-\frac{\epsilon^2j}{2}}O(j^{-\frac{1}{2}}),\\
\sqrt{j+1}-\sqrt{j}&=&\frac{1}{2}j^{-\frac{1}{2}}+O(j^{-\frac{3}{2}}).
\end{eqnarray*}
From the above argument, we have
\begin{eqnarray}\label{eq:3.9}
I_{31}=\frac{p\epsilon^p}{2-p}\sum_{j=1}^{\infty}\Big[\Phi(\epsilon\sqrt{j})
+\epsilon e^{-\frac{\epsilon^2j}{2}}O(j^{-\frac{1}{2}})
+\Phi(\epsilon\sqrt{j})O(j^{-1})+O(j^{-\frac{3}{2}})\Big].
\end{eqnarray}
Letting $y=\epsilon\sqrt{t}$, we have
\beqlb\label{d-3}
\frac{p\epsilon^p}{2-p}\sum_{j=1}^{\infty}\Phi(\epsilon\sqrt{j})=&&\frac{p\epsilon^p}{2-p}\int_{1}^{\infty}\Phi(\epsilon\sqrt{t})dt+O(\epsilon^p)\nonumber\\
&=&\frac{p\epsilon^{p-2}}{2-p}\int_{\epsilon}^{\infty}\Phi(y)2ydy+O(\epsilon^p)\nonumber\\
&=&\frac{p\epsilon^{p-2}\E [N^2]}{2-p}-\frac{p\epsilon^{p-2}}{2-p}\int_{0}^{\epsilon}2t\Phi(y)dy+O(\epsilon^p)\nonumber\\
&=&\frac{p\epsilon^{p-2}\E [N^2]}{2-p}+O(\epsilon^p).
\eeqlb
Similar to \eqref{d-3},
we can obtain that
\beqlb\label{d-4}
\frac{p\epsilon^{p+1}}{2-p}\sum_{j=1}^{\infty}O(j^{-\frac{1}{2}})e^{-\frac{\epsilon^2j}{2}}
=O(\epsilon^p),
\eeqlb
\beqlb\label{d-5}
&&\frac{p\epsilon^{p}}{2-p}\sum_{j=1}^{\infty}O(j^{-1})\Phi(\epsilon\sqrt{j})=O(\epsilon^p\log \frac{1}{\epsilon}),
\eeqlb
and
\beqlb\label{d-6}
&&\frac{p\epsilon^{p}}{2-p}\sum_{j=1}^{\infty}O(j^{-\frac{3}{2}})
=O(\epsilon^p).
\eeqlb

From \eqref{d-3} to \eqref{d-6}, we have
\begin{eqnarray}\label{eq:3.12}
I_{31}=\frac{p\epsilon^{p-2}\E [N^2]}{2-p}+O(\epsilon^p\log\frac{1}{\epsilon}).
\end{eqnarray}

For $I_{32}$, we have
\begin{eqnarray}\label{eq:3.14}
I_{32}&=&B_{-\frac{p}{2}}\int_{\epsilon}^{\infty}pt^{p-1}\Phi(t)dt
\nonumber\\
&=&B_{-\frac{p}{2}}\E [|N|^p]-B_{-\frac{p}{2}}\int_{0}^{\epsilon}pt^{p-1}\Phi(t)dt\nonumber
\\
&=&B_{-\frac{p}{2}}\E [|N|^p]+O(\epsilon^p).
\end{eqnarray}

Next, we deal with $I_{33}$.
Noting that $$(j+1)^{\frac{p}{2}}-j^{\frac{p}{2}}=O(j^{\frac{p}{2}-1}),\;\textrm{as}\;j\to\infty,$$
we have
\beqlb\label{d-7}\int_{\epsilon\sqrt{j}}^{\epsilon\sqrt{j+1}}pt^{p-1}\Phi(t)dt
\leq \epsilon^p\Phi(\epsilon\sqrt{j})[(j+1)^{\frac{p}{2}}-j^{\frac{p}{2}}]=\epsilon^p\Phi(\epsilon\sqrt{j})O(j^{\frac{p}{2}-1}),\eeqlb
since  $\Phi(t)$ is a decreasing function.

It follows from \eqref{d-7} that
\begin{eqnarray}
I_{33}&\leq&C\epsilon^p\sum_{j=1}^{\infty}j^{-1}\Phi(\epsilon\sqrt{j})\nonumber\\
&=&C\epsilon^p\int_{1}^{\infty}\frac{\Phi(\epsilon\sqrt{x})}{x}dx+O(\epsilon^p)
=C\epsilon^p\int_{\epsilon}^{\infty}\frac{\Phi(t)}{t}dt+O(\epsilon^p)\label{eq:3.15}\\
&\leq&C\epsilon^p\int_{\epsilon}^1t^{-1}dt+C\epsilon^p\int_{1}^{\infty}\Phi(t)dt+O(\epsilon^p)
=O(\epsilon^p\log \frac{1}{\epsilon}),\label{eq:3.16}
\end{eqnarray}
where the second equation in \eqref{eq:3.15} follows from the change of variable $t=\epsilon\sqrt{x}$ again.
\eqref{eq:3.12}, \eqref{eq:3.14} and \eqref{eq:3.16} imply that
\begin{eqnarray}\label{eq:3.18}
\sum_{n=1}^{\infty}\frac{1}{n^p}\int_{\epsilon n}^{\infty}px^{p-1}\P(|S_ n|\geq x)dx=
\frac{p\epsilon^{p-2}\E [N^2]}{2-p}+B_{-\frac{p}{2}}\E[|N|^p]+O(\epsilon^p\log\frac{1}{\epsilon}).
\end{eqnarray}
From \eqref{eq:3.5}, \eqref{eq:3.6} and \eqref{eq:3.18},
we get the part (a) in Proposition \ref{prop2.1}.

Below, we prove the part (b).   The proof of this part is similar to that of the part (a). However, some
modifications are needed to characterize the lower bound $\epsilon\sqrt{n\log n}$.
Note that, for any $\delta\in(0,1]$,
\begin{eqnarray}\label{eq:3.19}
\lambda_2(\epsilon,\delta)=&&\epsilon^2\sum_{n=2}^{\infty}\frac{(\log n)^{\delta}}{n}
\P(|S_n|\geq \epsilon\sqrt{n \log n})+\nonumber\\
&&\sum_{n=2}^{\infty}\frac{(\log n)^{\delta-1}}{n^2}
\int_{\epsilon\sqrt{n\log n}}^{\infty}2x\P(|S_n|\geq x)dx.
\end{eqnarray}
It follows from Kong \cite{K2015}  that
\begin{eqnarray}\label{eq:3.20}
\sum_{n=2}^{\infty}\frac{(\log n)^{\delta}}{n}
\P(|S_n|\geq \epsilon\sqrt{n \log n})=
\frac{\epsilon^{-2\delta-2}\E[|N|^{2\delta+2}]}{\delta+1}+O(1).
\end{eqnarray}
Hence, we only need to compute the second term in \eqref{eq:3.19}.
In fact, by the Fubini's theorem and
Lemma \ref{lem-3.3}, we have
\begin{eqnarray}
&&\sum_{n=2}^{\infty}\frac{(\log n)^{\delta-1}}{n^2}
\int_{\epsilon\sqrt{n\log n}}^{\infty}2x\P(|S_n|\geq x)dx\nonumber\\
&=&\sum_{j=2}^{\infty}\big(\sum_{n=2}^{j}\frac{(\log n)^{\delta-1}}{n}\big)
\int_{\epsilon\sqrt{\log j}}^{\epsilon\sqrt{\log (j+1)}}2t\Phi(t)dt\nonumber\\
&=:&I_{34}+I_{35}+I_{36},\label{eq:3.21}
\end{eqnarray}
where
\beqnn
I_{34}&&=\frac{1}{\delta}\sum_{j=2}^{\infty}(\log j)^{\delta}
\int_{\epsilon\sqrt{\log j}}^{\epsilon\sqrt{\log (j+1)}}2t\Phi(t)dt,
\\I_{35}&&=C_{\delta}\sum_{j=2}^{\infty}\int_{\epsilon\sqrt{\log j}}^{\epsilon\sqrt{\log (j+1)}}2t\Phi(t)dt,
\\I_{36}&&=\sum_{j=2}^{\infty}O(\frac{(\log j)^{\delta-1}}{j})
\int_{\epsilon\sqrt{\log j}}^{\epsilon\sqrt{\log (j+1)}}2t\Phi(t))dt.
\eeqnn

We first consider $I_{34}$. The integer mean theorem shows that
\beqlb\label{d-8}
\int_{\epsilon\sqrt{\log j}}^{\epsilon\sqrt{\log (j+1)}}2t\Phi(t)dt
=\epsilon^2\Phi(\epsilon\sqrt{\log\xi})(\log(j+1)-\log j)\;\textrm{for some}\; \xi\in(j,j+1).
\eeqlb
On the other hand, it follows from the Taylor expansion that
\begin{eqnarray*}
\Phi(\epsilon\sqrt{\log \xi})&=&\Phi(\epsilon\sqrt{\log j})+\epsilon O(j^{-1-\frac{\epsilon^2}{2}}(\log j)^{-\frac{1}{2}}),\\
\log(j+1)-\log j&=&j^{-1}+O(j^{-2}).
\end{eqnarray*}
By the above argument, we have
\begin{eqnarray*}
I_{34}&=&\frac{\epsilon^2}{\delta}\sum_{j=2}^{\infty}\frac{(\log j)^\delta}{j}
\Phi(\epsilon\sqrt{\log j})+O(\epsilon^2).
\end{eqnarray*}
Putting $t=\epsilon \sqrt{\log x}$, we have
\begin{eqnarray}\label{eq:3.22}
&&\frac{\epsilon^2}{\delta}\int_{2}^{\infty}\frac{(\log x)^\delta}{x}
\Phi(\epsilon\sqrt{\log x})dx\nonumber\\
&=&\frac{2\epsilon^{-2\delta}}{\delta}\int_{\epsilon\sqrt{\log 2}}^{\infty}t^{2\delta+1}
\Phi(t)dt=\frac{\epsilon^{-2\delta}\E[|N|^{2\delta+2}]}{\delta(\delta+1)}
+O(\epsilon^2),
\end{eqnarray}
since
 $$\int_{0}^{\epsilon\sqrt{\log 2}}2t^{2\delta+1}\Phi(t)dt\leq
\int_{0}^{\epsilon\sqrt{\log 2}}2t^{2\delta+1}dt=O(\epsilon^{2\delta+2}).$$
It follows from \eqref{eq:3.22} that
\begin{eqnarray}\label{eq:3.22a}
I_{34}\leq\frac{\epsilon^{-2\delta}
\E[|N|^{2\delta+2}]}{\delta(\delta+1)}
+O(\epsilon^2).
\end{eqnarray}

For $I_{35}$, we have
\begin{eqnarray}\label{eq:3.23}
I_{35}=C_{\delta}\int_{0}^{\infty}2t\Phi(t)dt
-C_{\delta}\int_{0}^{\epsilon\sqrt{\log 2}}2t\Phi(t)dt
=C_{\delta}\E[N^2]+O(\epsilon^2).
\end{eqnarray}

Finally, we look at $I_{36}$.  By \eqref{d-8}, we have
\begin{eqnarray}\label{eq:3.24}
I_{36}\leq C \epsilon^2 \sum_{j=2}^{\infty}\frac{(\log j)^{\delta-1}}{j^2}
\Phi(\epsilon\sqrt{\log j})+O(\epsilon^2)=O(\epsilon^2),
\end{eqnarray}
where
$$0<\delta\leq 1, \mbox{ and }\Phi(\epsilon\sqrt{\log j})\leq 1 \mbox{ for any }j\geq 2.$$

Hence, from \eqref{eq:3.22a} to \eqref{eq:3.24},  we obtain
\begin{eqnarray}\label{eq:c}
&&\sum_{n=2}^{\infty}\frac{(\log n)^{\delta-1}}{n^2}
\int_{\epsilon\sqrt{n\log n}}^{\infty}2x\P(|S_n|\geq x)dx\nonumber
\\&& \qquad=\frac{\epsilon^{-2\delta}
\E[|N|^{2\delta+2}]}{\delta(\delta+1)}+C_{\delta}\E[N^2]+O(\epsilon^2).
\end{eqnarray}
Therefore, we get from \eqref{eq:3.19}, \eqref{eq:3.20} and \eqref{eq:c} that
the  part (b)  holds.  The proof of this proposition is finished.

\qed

Next, we focus on the proof of Proposition \ref{prop2.2}.

\noindent{\it Proof of Proposition \ref{prop2.2}:}
Let $\{X,X_n;~n\geq 1\}$ be a sequence of i.i.d. random
variables with mean zero, $\E [X^2]=1$
and $\E[|X|^q]<\infty$ for some $q\in (2,3]$.  Moreover, for any $n\in \N$, set
\begin{eqnarray}\label{eq:3.25}
\Delta_n:=\sup_{t}\Big|\P(|S_n|\geq \sqrt{n}t)-\Phi(t)\Big|.
\end{eqnarray}

We first prove the part (a).
 For any $M\geq 1$,
let $$H_1(\epsilon)=M\epsilon^{-2\gamma},$$
where $\gamma=\frac{q-p}{2q-2-p}>0$ with $0<p<2$ and $2<q\leq 3$.

By the change of variable $t=\frac{x}{\sqrt{n}}$, we have
\begin{eqnarray}\label{eq:3.30}
&&\sum_{n=1}^{\infty}n^{-p}\int_{\epsilon n}^{\infty}px^{p-1}\Big|
\P(|S_n|\geq x)-\Phi(\frac{x}{\sqrt{n}})\Big|dx\nonumber\\
&&\qquad=\sum_{n=1}^{\infty}n^{-\frac{p}{2}}\int_{\epsilon \sqrt{n}}^{\infty}pt^{p-1}\Big|
\P(|S_n|\geq \sqrt{n}t)-\Phi(t)\Big|dt.
\end{eqnarray}
Based on $H_1(\epsilon)$, we split the right hand side of \eqref{eq:3.30}  into two parts.
We first observe the first part, that is,
\begin{eqnarray}\label{eq:3.31}
\sum_{n\leq H_1(\epsilon)}n^{-\frac{p}{2}}\int_{\epsilon \sqrt{n}}^{\infty}pt^{p-1}\Big|
\P(|S_n|\geq \sqrt{n}t)-\Phi(t)\Big|dt\leq \sum_{n\leq H_1(\epsilon)}n^{-\frac{p}{2}}(J_1+J_2),
\end{eqnarray}
where
$$
J_1=\int_{0}^{\Delta_n^{-\frac{1}{2p}}}pt^{p-1}\Big|\P(|S_n|\geq \sqrt{n}t)-\Phi(t)\Big|dt
$$
and
$$
J_2=\int_{\Delta_n^{-\frac{1}{2p}}}^{\infty}pt^{p-1}\Big|\P(|S_n|\geq \sqrt{n}t)-\Phi(t)\Big|dt.
$$
It follows from $$\frac{S_n}{\sqrt{n}}\overset{d}{\rightarrow}N$$ that
$$\Delta_n\rightarrow0$$as $n\rightarrow\infty$. Thus,
\beqlb\label{d-10}J_1\leq\Delta_n\int_{0}^{\Delta_n^{-\frac{1}{2p}}}pt^{p-1}dt
=\Delta_n^{\frac{1}{2}}\rightarrow 0~~(n\rightarrow \infty).\eeqlb

On the other hand, it follows from $\frac{S_n}{\sqrt{n}}\overset{d}{\rightarrow}N$
and Lemma \ref{lem-3.1} that for large enough
$t$,
\beqlb\label{d-9}
\Big|\P(|S_n|\geq \sqrt{n}t)-\Phi(t)\Big|\leq Ce^{-\frac{t^2}{2}}t^{-1}\leq Ce^{-t}t^{1-p}.
\eeqlb
Hence,  as $n$ goes to $\infty$,
\beqlb\label{d-11}J_2\leq C\int_{\Delta_n^{-\frac{1}{2p}}}^{\infty}e^{-t}dt\rightarrow0.\eeqlb
By the Toeplitz lemma, \eqref{d-10} and \eqref{d-11}, we have that for any $M\geq 1$,
\begin{eqnarray}\label{eq:3.32}
\lim_{\epsilon\searrow 0}\epsilon^{\gamma(2-p)}\sum_{n\leq H_1(\epsilon)}n^{-\frac{p}{2}}(J_1+J_2)
=\lim_{\epsilon\searrow 0}
\frac{M^{1-\frac{p}{2}}}{(H_1(\epsilon))^{1-\frac{p}{2}}}\sum_{n\leq H_1(\epsilon)}n^{-\frac{p}{2}}(J_1+J_2)=0.
\end{eqnarray}

Next, we look at the second part.  Recall that  Bikjalis \cite{B1966} got the following non-uniform large deviation estimate
, for any $x>0$,
\begin{eqnarray}\label{eq:3.29}
\Big|\P(S_n>\sqrt{n}x)-\P(N>x)\Big|
\leq\frac{C\E[|X|^q]}{n^{\frac{q}{2}-1}(1+x^q)},~~2<q\leq 3.
\end{eqnarray}
Hence
\begin{eqnarray*}
&&\sum_{n\geq H_1(\epsilon)}n^{-\frac{p}{2}}\int_{\epsilon \sqrt{n}}^{\infty}pt^{p-1}\Big|
\P(|S_n|\geq \sqrt{n}t)-\Phi (t)\Big|dt\nonumber\\
&&\qquad\qquad\leq\sum_{n\geq H_1(\epsilon)}n^{-\frac{p}{2}}\int_{\epsilon \sqrt{n}}^{\infty}
\frac{Ct^{p-1}}{n^{\frac{q}{2}-1}(1+t^q)}dt\nonumber\\
&&\qquad\qquad\leq\sum_{n\geq H_1(\epsilon)}C n^{1-\frac{p}{2}-\frac{q}{2}}\int_{\epsilon \sqrt{n}}^{\infty}
t^{p-q-1}dt=C\epsilon^{p-q}
\sum_{n\geq H_1(\epsilon)}n^{1-q}.
\end{eqnarray*}
Thus,
\begin{eqnarray}\label{eq:3.33}
&&\limsup_{\epsilon\searrow0}\epsilon^{\gamma(2-p)}
\sum_{n\geq H_1(\epsilon)}n^{-\frac{p}{2}}\int_{\epsilon \sqrt{n}}^{\infty}pt^{p-1}\Big|
\P(|S_n|\geq \sqrt{n}t)-\Phi (t)\Big|dt\nonumber\\
&&\qquad\qquad \leq C\limsup_{\epsilon\searrow0}\epsilon^{\gamma(2-p)+(p-q)}
 (H_1(\epsilon))^{2-q}\nonumber
 \\
 &&\qquad\qquad=CM^{2-q}\searrow0,\;\mbox{ as }\;M\nearrow\infty,
\end{eqnarray}
since $q>2$.

 Finally,  we get from \eqref{eq:3.32} and \eqref{eq:3.33} that the part (a).

Below, we prove the part (b).  In fact,
the proof of this part is similar to that of part (a)
and we only write the
modifications in the following.

For any $M\geq 1$, define $$H_2(\epsilon)=M\epsilon^{-2}(\log\frac{1}{\epsilon})^{\frac{2\delta}{2-q}}.$$
It is obvious that as $\epsilon\searrow0$, $$H_2(\epsilon)\nearrow\infty.$$

The change of variable $t=\frac{x}{\sqrt{n}}$ yields that
\begin{eqnarray}\label{eq:3.34}
&&\sum_{n=2}^{\infty}\frac{(\log n)^{\delta-1}}{n^2}\int_{\epsilon\sqrt{n\log n}}^{\infty}
2x\Big|\P(|S_n|\geq x)-\Phi\big(\frac{x}{\sqrt{n}}\big)\Big|dx\nonumber\\
&&\qquad\qquad=\sum_{n=2}^{\infty}\frac{(\log n)^{\delta-1}}{n}\int_{\epsilon\sqrt{\log n}}^{\infty}
2t\Big|\P(|S_n|\geq \sqrt{n}t )-\Phi(t)\Big|dt.
\end{eqnarray}

According to $H_2(\epsilon)$, we split the sum in \eqref{eq:3.34}  into two part. For the first part,
\begin{eqnarray*}
\sum_{n=2}^{H_2(\epsilon)}\frac{(\log n)^{\delta-1}}{n}\int_{\epsilon\sqrt{\log n}}^{\infty}
2t\Big|\P(|S_n|\geq \sqrt{n}t )-\Phi(t)\Big|dt
\leq\sum_{n=2}^{H_2(\epsilon)}\frac{(\log n)^{\delta-1}}{n}(J_3+J_4),
\end{eqnarray*}
where
$$
J_3=\int_{0}^{\Delta_n^{-\frac{1}{4}}}2t\Big|\P(|S_n|\geq \sqrt{n}t)-\Phi(t)\Big|dt
$$
and
$$
J_4=\int_{\Delta_n^{-\frac{1}{4}}}^{\infty}2t\Big|\P(|S_n|\geq \sqrt{n}t)-\Phi(t)\Big|dt
$$
with $\Delta_n$ being defined by \eqref{eq:3.25}.
Similar to \eqref{d-10} and \eqref{d-11}, we  obtain that $J_3$ and $ J_4$ goes to $ 0$
as $n\rightarrow \infty$, respectively.  Thus, by using the Toeplitz lemma, we have, for any  $M\geq 1$,
\begin{eqnarray}\label{eq:3.35}
&&\lim_{\epsilon\searrow0}(\log\frac{1}{\epsilon})^{-\delta}
\sum_{n=2}^{H_2(\epsilon)}\frac{(\log n)^{\delta-1}}{n}(J_3+J_4)\nonumber\\
&&\qquad\qquad=\lim_{\epsilon\searrow0}\frac{2^{\delta}}{[\log H_2(\epsilon)]^{\delta}}
\sum_{n=2}^{H_2(\epsilon)}\frac{(\log n)^{\delta-1}}{n}(J_3+J_4)=0.
\end{eqnarray}
Furthermore, by using the large deviation \eqref{eq:3.29}, we have
\begin{eqnarray}
&&\sum_{n\geq H_2(\epsilon)}\frac{(\log n)^{\delta-1}}{n}\int_{\epsilon\sqrt{\log n}}^{\infty}
2t\Big|\P(|S_n|\geq \sqrt{n}t )-\Phi(t)\Big|dt\nonumber\\
&&\qquad\qquad\leq\sum_{n\geq H_2(\epsilon)}\frac{(\log n)^{\delta-1}}
{n^{\frac{q}{2}}}\int_{\epsilon\sqrt{\log n}}^{\infty}
\frac{C t\E[|X|^q]}{1+t^q}dt\nonumber\\
&&\qquad\qquad=C\epsilon^{2-q}
\sum_{n\geq H_2(\epsilon)}(\log n)^{\delta-\frac{q}{2}}
n^{-\frac{q}{2}}\label{eq:3.36}.
\end{eqnarray}
Note that $$0<\delta\leq 1\mbox{ and }2<q\leq3.$$ Hence, we have
$\delta-\frac{q}{2}<0$
and
\begin{eqnarray*}
&&\limsup_{\epsilon\searrow0}(\log\frac{1}{\epsilon})^{-\delta}
\sum_{n\geq H_2(\epsilon)}\frac{(\log n)^{\delta-1}}{n}\int_{\epsilon\sqrt{\log n}}^{\infty}
2t\Big|\P(|S_n|\geq \sqrt{n}t )-\Phi(t)\Big|dt\\
&&\qquad\qquad\leq C\limsup_{\epsilon\searrow0}(\log\frac{1}{\epsilon})^{-\delta}\epsilon^{2-q}
\sum_{n\geq H_2(\epsilon)}n^{-\frac{q}{2}}\nonumber\\
&&\qquad\qquad \leq C\limsup_{\epsilon\searrow0}(\log\frac{1}{\epsilon})^{-\delta}\epsilon^{2-q}
 (H_2(\epsilon))^{1-\frac{q}{2}}\nonumber
 \\&&\qquad\qquad \leq CM^{1-\frac{q}{2}}\searrow0,\;\mbox{as }M\nearrow\infty.
\end{eqnarray*}
Thus, we get from \eqref{eq:3.35} and  \eqref{eq:3.35} the part (b) of Proposition \ref{prop2.2}.
\qed

Finally, we prove the main result of this paper.

\noindent{\it Proof of Theorem 2.1:}  In order to prove it, we first should point out that \eqref{eq:3.5} and \eqref{eq:3.19} also hold under the conditions of Theorem \ref{thm-2}.  We prove this theorem by two steps.

We first prove the part (a) of this theorem.  Recall that
He and Xie \cite{HX2013}  obtained that
\beqlb\label{d-12}
\sum_{n=1}^{\infty}\P(|S_n|\geq \epsilon n)=\epsilon^{-2}(\sigma^2+o(\epsilon^{q-2})).
\eeqlb
Combing \eqref{eq:3.5} and \eqref{d-12}, we get that  for any $0<p<2$
\begin{eqnarray*}
\lambda_{1}(\epsilon,p)=\sigma^2\epsilon^{p-2}+\epsilon^{p+q-4}o(1)+
\sum_{n=1}^{\infty}n^{-p}\int_{\epsilon n}^{\infty}px^{p-1}
\P(|S_n|\geq x)dx.
\end{eqnarray*}
Thus, by \eqref{eq:3.18}, we have
\begin{eqnarray}\label{eq:a}
\lambda_{1}(\epsilon,p)&=&\frac{2\sigma^2\epsilon^{p-2}}{2-p}+B_{-\frac{p}{2}}\sigma^p\E[|N|^p]
+\epsilon^{p+q-4}o(1)+O(\epsilon^p\log\frac{1}{\epsilon})\nonumber\\
&&~~~~~~~~~~~+\sum_{n=1}^{\infty}n^{-p}\int_{\epsilon n}^{\infty}px^{p-1}
\Big[\P(|S_n|\geq x)-\Phi(\frac{x}{\sqrt{n}})\Big]dx.
\end{eqnarray}
Hence,
combing \eqref{eq:a} and part (a) of Proposition \ref{prop2.2} together,
we get that
$$\lim_{\epsilon\searrow0}\epsilon^{\gamma(2-p)}\Big[\lambda_1(\epsilon,p)-\frac{2\sigma^2\epsilon^{p-2}
}{2-p}\Big]=0,$$
where $$0<p<2,~~0<\gamma=\frac{q-p}{2q-2-p}<1,$$
and $$ \gamma(2-p)+p+q-4=\frac{2(q-2)^2}{2q-2-p}>0.$$
Hence the  part (a)  holds.

Next, we prove the part (b). The arguments for the part (b) are similar to those for the part (a).
Kong \cite{K2015} got that \eqref{eq:3.20} also holds
for all random variables satisfying \eqref{eq:2.2a}.
Thus, by \eqref{eq:c} and part (b) of Proposition \ref{prop2.2}, we  get that
$$\lim_{\epsilon\searrow0}(\log\frac{1}{\epsilon})^{-\delta}\Big[\lambda_{2}(\epsilon,\delta)
-\frac{\sigma^{2\delta+2}}{\delta}\E [|N|^{2\delta+2}]\epsilon^{-2\delta}\Big]=0.$$
 Therefore, we get part (b).  The proof of Theorem \ref{thm-2} is finished.
\qed

\noindent{\bf Acknowledgments:}\  This work was  supported by the Natural Science Foundation of China (No.11361007), the Guangxi Natural Science Foundation (Nos.2012GXNSFBA053010 and 2014GXNSFCA118001) and  the Project for Fostering Distinguished Youth Scholars of Shandong University of Finance and Economics.

\end{document}